\documentclass{amsart}%
\usepackage{graphicx}
\usepackage{amsmath}
\usepackage{amscd}
\usepackage{amsfonts}

\newtheorem{theorem}{Theorem}[section]
\newtheorem{lemma}[theorem]{Lemma}

\newtheorem{example}[theorem]{Example}

\newtheorem{corollary}[theorem]{Corollary}
\newtheorem{remark}[theorem]{Remark}

\newcommand{\be}{\begin{equation}}
\newcommand{\ee}{\end{equation}}
\newcommand{\bes}{\begin{equation*}}
\newcommand{\ees}{\end{equation*}}

\newcommand{\cX}{\mathcal{X}}

\newcommand{\cS}{\mathcal{S}}

\newcommand{\Rp}{\mathbb{R}_+}

\begin{document}

\title{Continuous extension of a densely parameterized semigroup}

\author{Eliahu Levy}
\address{Department of Mathematics, Technion - Israel
Institute of Technology, 32000, Haifa, Israel.}
\email{eliahu@techunix.technion.ac.il}
\author{Orr Moshe Shalit}
\address{Department of Mathematics, Technion - Israel
Institute of Technology, 32000, Haifa, Israel.}
\email{orrms@techunix.technion.ac.il}
\thanks{O.M.S. was partially
supported by the Gutwirth Fellowship.}
\keywords{Contraction semigroup, nonlinear semigroup
extension, densely defined.}
\subjclass[2000]{47D03, 47H20.}
%\author{}
\date{January 30, 2007}

\begin{abstract}
{Let $\cS$ be a dense sub-semigroup of $\Rp$, 
and let $X$ be a separable, reflexive Banach space. 
This note contains a proof that every weakly continuous contractive semigroup 
of operators on $X$ over $\cS$ can be extended to a weakly continuous semigroup over $\Rp$. 
We obtain similar results for nonlinear, nonexpansive semigroups as well. 
As a corollary we characterize all densely parametrized semigroups which are extendable to semigroups over $\Rp$.
}
\end{abstract}
\maketitle

\section{Introduction}
Let $X$ be a Banach space, and let $\cS$ be a dense sub-semigroup of
$\Rp$. A \emph{semigroup of operators over $\cS$} is a family $T =
\{T_s\}_{s\in\cS}$ of operators on $X$ such that
$$T_{s+t} = T_s \circ T_t \,\, , \,\, s,t \in \cS .$$
If $0 \in \cS$, we also require that $T_0 = I$. 
We shall refer below to such a semigroup 
as a \emph{densely parametrized} semigroup. 

The word
\emph{operator} shall mean henceforth \emph{linear operator}
unless otherwise stated. 
A semigroup $T$ (over $\cS$)
is said to be \emph{weakly continuous} if for all $x\in X, y \in
X^*$, the function $\cS\ni s \mapsto y(T_s(x))$ is a continuous
function. Left and right weak continuity are defined similarly.

The theory of weakly continuous semigroups over $\Rp$ is highly developed \cite{HP}. Some of
the techniques used for semigroups over $\Rp$ cannot be used when
one considers a semigroup of operators over an arbitrary semigroup
$\cS$. For example, the existence of a generator for the semigroup
can be proved using Bochner integration. But if one has a
semigroup of operators, say, over the rationals, then one cannot
integrate. The main result of this paper is that if $\cS$ is a
dense sub-semigroup of $\Rp$ and $X$ is a separable, reflexive
Banach space, then every right weakly continuous contractive
semigroup on $X$ over $\cS$ can be extended to a weakly continuous
semigroup over $\Rp$. 

A similar but weaker result is also obtained for 
semigroups of nonlinear operators. A nonlinear map $A$ is said to be \emph{nonexpansive}
if $A$ is Lipschitz continuous with a Lipschitz constant, denoted by $\|A\|$, that is not larger than $1$.
We shall show that, under the same assumptions on $X$ and $\cS$, every right weakly continuous
semigroup of nonexpansive maps that are continuous with respect to the weak topology on $X$ can be extended to a right weakly continuous semigroup over $\Rp$. 

The result that every densely parametrized semigroup of (linear) contractions
that is weakly continuous from the right may be extended to a continuous semigroup
parametrized by $\Rp$ may seem rather expected. Indeed, if the semigroup is assumed to 
be \emph{strongly} continuous from the right, that is, if for all $x\in X$ the function $\cS\ni s \mapsto T_s(x)$ is continuous from the right (where $X$ is given the norm topology), 
then constructing a continuous extension is straightforward. One is tempted to think that a densely parametrized
semigroup that is continuous with respect to any reasonable topology can always be extended to a continuous
semigroup (with respect to the same topology) over $\Rp$. The following example may serve to illustrate
that things do not always work as expected.
\begin{example}
\emph{
Let $X$ be the closed subspace of $L^\infty(\mathbb{R})$ spanned by the functions $x \mapsto e^{iqx}$ with $q \in \mathbb{Q}$. We endow $X$ with the topology inherited from the weak-$*$ topology on $L^\infty(\mathbb{R})$.
We call this topology the $L^1$ weak topology on $X$. Let $T = \{T_s\}_{s \in \mathbb{Q}}$ be a group of isometric multiplication operators on $X$ given by 
\bes
(T_s f)(x) = e^{isx}f(x).
\ees
For every $f \in X$, the function $s \mapsto T_s f$ is continuous with respect to the $L^1$ weak topology,
but $T$ cannot be extended to an $L^1$ weakly continuous semigroup over $\mathbb{R}$. Indeed,
if $T$ was extendable then for $r \notin \mathbb{Q}$ and for all $g \in L^1$ we would have
\bes
\lim_{s \rightarrow r} \int_{\mathbb{R}} g(x)e^{isx}f(x)dx = \int_{\mathbb{R}} g(x)(T_r f)(x)dx ,
\ees
from which it follows (using Lebesgue's Dominated Convergence Theorem) that $T_r$ must be given by multiplication by
$e^{irx}$. However, $X$ is not closed under multiplication by $e^{irx}$.}
\end{example}

\section{The main result}
Throughout this section, let $X$ be a separable and reflexive
Banach space, with a dual $X^*$, and let $\cS$ be a dense
sub-semigroup of $\Rp = [0,\infty)$. A \emph{contractive semigroup
on $X$} (over $\cS$) is simply a semigroup $T = \{T_s\}_{s\in\cS}$
such that $T_s$ is a contraction for all $s\in\cS$, that is, $\|T_s\|$ is a linear operator
such that $\|T_s\| \leq 1$.
%A \emph{contractive semigroup
%on $X$} (over $\cS$) is a family $T = \{T_s\}_{s\in\cS}$ of
%contractions satisfying $T_{s+t} = T_s \circ T_t$ for all
%$s,t\in\cS$. If $0 \in \cS$, we also require that $T_0 = I$. $T$
%is said to be \emph{weakly continuous} if for all $x\in X, y\in
%X^*$, the function $s \mapsto y(T_s (x))$ is continuous.

Recall that the pair $(X,X^*)$ satisfies:
\begin{equation}\label{eq:norm}
\|x\| = \max_{y \in X^*, \|y\|=1} |y(x)|,
\end{equation}
and that $X$ is \emph{weakly sequentially complete}, that is,
it has the property: if $\{x_n\} \subset X$ is such that
for all $y\in X^*$, $\{y(x_n)\}$ converges, then there is $x\in X$ such that $y(x_n) \rightarrow y(x)$ for all
$y \in X^*$.

\begin{theorem}\label{thm:weak}
Let $X$ and $\cS$ be as above. Let $T = \{T_s\}_{s\in\cS}$ be a
contractive semigroup on $X$, such that
\begin{equation}\label{eq:continuity}
\lim_{\cS\ni s\rightarrow 0^+} y(T_s(x)) = y(x) \,\, , \,\, x \in X, y \in X^* .
\end{equation}
Then $T$ can be extended to a weakly continuous contractive
semigroup $\{T_t\}_{t\geq 0}$.
\end{theorem}
\begin{remark}
\emph{If $T$ is a \emph{nonlinear} semigroup of nonexpansive
maps satisfying, in addition to the above conditions, the
assumption that for all $s\in \cS$, $T_s$ is continuous in the
weak topology of $X$, then the following proof will guarantee that
we can extend $T$ to a \emph{right} weakly continuous semigroup over $\Rp$ of
nonexpansive maps. Throughout the proof, we shall indicate where the differences between
linear and nonlinear semigroups  occur.}
\end{remark}
\begin{proof} We shall split the proof into a number of logical steps.

\vspace{0.3cm}
\noindent{\bf 1. Simplifying assumptions.}

We assume that $X$ is a real Banach space,
as the complex case follows easily by considering
the real and imaginary parts of the functionals appearing in the proof.
%We may also assume that $\cS$ is countable, otherwise we look at a countable
%subsemigroup $\cS'$. This is possible by taking a sequence $\{s_n\}$ converging to $0$,
%and letting $\cS$ be the subsemigroup this sequence generates. Since $T$ is continuous on
%$\cS$ and on $\cS'$, the extension of $T$ from $\cS'$ to $\Rp$ will agree with
%the original definition of $T$ on $\cS$.
%
%Finally, we assume that $T$ is right continuous at any $s\in \cS$, as this clearly follows from (\ref{eq:continuity}).
We also assume that $T$ is right continuous at any $s\in \cS$, as
this clearly follows from (\ref{eq:continuity}).

\vspace{0.3cm}
\noindent{\bf 2. Preliminary definitions.}

For any (real valued) continuous function $\varphi$ on $\cS$ we define a function $\varphi^-$ on $\Rp$ by
\bes
\varphi^- (t) = \inf\{h(t) : h\in RUSC(\Rp), \forall s\in\cS . \varphi(s) \leq h(s)\}
\ees
for all $t \in \Rp$, where $RUSC(\Rp)$ denotes the space of right upper-semicontinuous%
\footnote{A \emph{right upper semicontinuous} function is just an
upper semicontinuous function with respect to the half-open topology
generated on $\Rp$ by the half open intervals of the type: $[a,b)$.

Note, that the open sets for the latter topology are characterized as those
whose connected components (with respect to the usual topology) are all intervals
open above, necessarily at most countable in number. Thus any set open for the
half-open topology turns into a usual open set by deleting an at most contable
set of points, hence the half-open interior (resp.\ closure) of any set differs
from the usual one by an at most countable set. One concludes that the properties
of a set being dense, resp.\ Baire, meager, residual, coincide for the half-open
and the usual topologies. In particular, $\Rp$ with the half-open topology is
a Baire space.}
(RUSC) functions on $\Rp$.
Similarly, we define $\varphi_-$ as the supremum of all right lower-semicontinuous functions (RLSC) smaller than
$\varphi$. It is clear that $\varphi_- \leq \varphi \leq \varphi^-$, $\varphi^-$ is RUSC, and $\varphi_-$ is RLSC.

For every fixed $x \in X, y \in X^*$ we can define a right continuous function on $\cS$ by
\be\label{eq:deff}
f(s;x,y) = y(T_s (x)) .
\ee
Our aim is to prove
\be\label{eq:what_to_show}
\left(f(t;x,y)\right)^- = \left(f(t;x,y)\right)_- \,\, , \,\, t\in \Rp, x\in X,y\in X^*.
\ee
Before we do that, we concentrate in the next two steps to show how the theorem follows from this fact.

\vspace{0.3cm}
\noindent{\bf 3. Showing how (\ref{eq:what_to_show}) gives rise to a weakly right-continuous contractive semigroup.}

Define
\be\label{eq:E}
E = \{t\in\Rp : \forall x\in X,y\in X^* . \left(f(t;x,y)\right)^- = \left(f(t;x,y)\right)_- \} .
\ee

Observe that $\cS
\subseteq E$. This follows from the fact that for all fixed $s\in
\cS$, $x\in X$ and $y\in X^*$, the functions
$$(y(T_s(x))+\epsilon)\cdot \chi_{[s,s+\delta)} + \infty\cdot\chi_{[s,s+\delta)^c}$$
and
$$(y(T_s(x))-\epsilon)\cdot \chi_{[s,s+\delta)} - \infty\cdot\chi_{[s,s+\delta)^c}$$
are right continuous, and for some $\delta>0$ they dominate and are dominated by the function  $\cS \ni t \mapsto f(t;x,y)$, respectively. We then have $f(s;x,y)^- - f(s;x,y)_- < 2\epsilon$, for all $\epsilon$, so $s \in E$.

For any $t\in \Rp$, if $\cS\ni s_n \searrow t$, then for all $x\in X, y\in X^*$,
\begin{align*}
& \left(f(t;x,y)\right)_- \leq \liminf \left(f(s_n;x,y)\right)_- = \liminf y(T_{s_n}(x)) \leq\\
& \leq \limsup y(T_{s_n}(x)) = \limsup \left(f(s_n;x,y)\right)^- \leq \left(f(t;x,y)\right)^- ,
\end{align*}
because $\left(f(\cdot;x,y)\right)_-$ is RLSC and $\left(f(\cdot;x,y\right)^-$ is RUSC.
If $t\in E$, then this means that $y(T_{s_n}(x)) \longrightarrow \left(f(t;x,y)\right)^-$ regardless of the choice of
$\{s_n\}$ and for all $x\in X$ and $y\in X^*$. Thus for $t\in E \setminus \cS$ we may define $T_t x$ to be the
weak limit $\lim_n T_{s_n} x$, where $\{s_n\}$ is any sequence in $\cS$ converging to $t$ from the right (this is where we
use the fact that $X$ is weakly sequentially complete). Note that
for $t\in E \cap \cS$ this weak limit would turn out to be the same $T_t$ that we started with. We will use this
below before we shall actually prove that $E = \Rp$.

Now if $E = \Rp$, then we get a family $\{T_t\}_{t\geq 0}$ of
linear operators on $X$. Equation (\ref{eq:norm}) implies that the
operators in this family are contractions. $\{T_t\}_{t\geq 0}$ is
weakly continuous from the right, since $y(T_t(x)) =
\left(f(t;x,y) \right)^- = \left(f(t;x,y) \right)_-$, a
right-continuous function in $t$. Also, in either case $0\in \cS$
or $0\notin\cS$, $T_0 = I$ by assumption.

To show that $\{T_t\}_{t\geq 0}$ is a semigroup, we first show
that \be\label{eq:halfSG} T_{s+t} = T_s \circ T_t \,\, , \,\, s
\in \cS, t \in \Rp. \ee Let $\cS\ni t_n \searrow t$, and fix $x\in
X, y\in X^*$. On one hand \be\label{eq:T converges} y(T_s\circ
T_{t_n}(x)) = y(T_{s+t_n}(x)) \longrightarrow y(T_{s+t}(x)) . \ee
On the other hand, \bes y(T_s\circ T_{t_n}(x)) =
y(T_s(T_{t_n}(x))) \longrightarrow y(T_s(T_t(x))) = y(T_s \circ
T_t (x)) , \ees because $T_{t_n}(x)$ converges weakly to
$T_{t}(x)$, and $T_s$ is continuous in the weak topology (as any
bounded operator. This is the main reason why in the nonlinear
case we assume that $T_s$ is weakly continuous, for all $s\in
\cS$). Together with (\ref{eq:T converges}) and (\ref{eq:norm}),
this means that (\ref{eq:halfSG}) holds.

Now let $s,t\in \Rp$, and let $\cS\ni s_n \searrow s$. On one
hand, from equation (\ref{eq:halfSG}) and the weak right
continuity of $\{T_t\}_{t\geq 0}$, it follows that for all $x\in
X$ \bes y(T_{s_n} \circ T_{t} (x)) = y(T_{s_n + t}(x))
\longrightarrow y(T_{s+t} (x)). \ees On the other hand, for all
$x$ and $y$, \bes y(T_{s_n}(T_t(x))) \rightarrow y(T_{s}(T_t(x)))
, \ees where we used again the weak right-continuity of
$\{T_t\}_{t\geq 0}$. Thus \bes T_{s+t} = T_s \circ T_t \,\, , \,\,
s, t \in \Rp. \ees

\vspace{0.3cm} \noindent{\bf 4. $\{T_t\}_{t\geq 0}$ (once defined)
is two sided weakly continuous.}

From the previous step, it follows that the semigroup $T$ extends
to a right weakly continuous contractive semigroup which we shall
also call $T$. It follows from classical results that $T$ is
weakly (and, in fact, strongly) continuous from the left as well
(see the corollary on page 306, \cite{HP}. This step does not hold
for the nonlinear case).

\vspace{0.3cm} \noindent{\bf 5. Two Lemmas.}

In this step we prove two technical lemmas, in order to make the
main parts of the proof smoother.

\begin{lemma}\label{lem:subspace}
For every $t\in\Rp$, the set
$$A_t = \{(x,y) \in X \times X^* :  \left(f(t;x,y)\right)^- = \left(f(t;x,y)\right)_- \}$$
is closed in $X \times X^*$.
\end{lemma}
\begin{proof}
Let $(x,y) \in \overline{A_t}$. We shall show that
$f(t;x,y)^- - f(t;x,y)_- \leq \epsilon$ for every $\epsilon > 0$. Indeed, given $\epsilon \in (0,1)$, let $(w,z) \in A_t$ 
such that  $\|w-x\|, \|z-y\| < \frac{\epsilon}{6 (M + N)}$, where $M:= \max\{\|x\|,\|y\|\} + 1$, and
$N$ is a bound for $\|T_s(x)\|$ for all $s \in \cS\cap[0,t+1]$. The existence of such a bound $N$ follows from 
(\ref{eq:continuity}), together with the semigroup property and the Principle of Uniform Boundedness  
(of course, if $T$ is a semigroup of linear operators, $N$ can be taken to be $\|x\|$).

Because $(w,z)\in A_t$, there is a $\delta\in(0,1)$ such that for all $s\in [t,t+\delta)$,
$$f(t;w,z)_- - \epsilon/6 < f(s;w,z) < f(t;w,z)^- + \epsilon/6 .$$
But then for all $s \in [t,t+\delta)$
\begin{align*}
f(s;x,y) &= y(T_s(w)) + y(T_s(x)) - y(T_s(w))\\
&\leq y(T_s(w)) + \|y\|\|T_s\|\|x-w\| \\
&\leq z(T_s(w)) + (y-z)(T_s(w)) + \epsilon/6 \\
&< f(t;w,z)^- + \epsilon/2 .
\end{align*}
Similarly, for all such $s$, $f(s;x,y)
> f(t;w,z)_- - \epsilon/2$. It follows that $f(t;x,y)^- - f(t;x,y)_- \leq
\epsilon$, for all $\epsilon$, in other words, $(x,y)\in A_t$.
\end{proof}

\begin{lemma}\label{lem:varphipsi}
Let $\varphi,\psi: \cS \rightarrow \mathbb{R}$ be right continuous, and let $c \in \Rp$ be such that
$$\varphi^-(s+c) \leq \psi(s) \,\, , \,\, s\in \cS .$$
Then
\be\label{eq:varphi}
\varphi^-(t+c) \leq \psi^-(t) \,\, , \,\, t\in \Rp .
\ee
A similar statement, with inequalities reversed and using $\varphi_-,\psi_-$ instead of $\varphi^-,\psi^-$, is also true.
\end{lemma}
\begin{proof}
Let $h$ be a RUSC function dominating $\psi$ on $\cS$. Then the
function $h_c$ given by $h_c(t) = h(t-c)$ for $t\geq c$, and
$h_c(t) = \infty$ for $t<c$, is RUSC and dominates $\varphi^-$ on
$c+\cS$. Let $c\leq s \in \cS$, and let $c + s_n \in c+\cS$ such
that $c + s_n \searrow s$. Since 
$\varphi$ is right continuous at
$s$, we have
$$\varphi(s)  = \lim_{n} \varphi^-(s_n + c) \leq \lim\sup_{n} \psi(s_n) \leq \lim\sup_{n} h(s_n) \leq h_c(s) .$$
Thus, $h_c$ is RUSC and dominates $\varphi^-$ on $\cS$, so
$\varphi^-(t+c) \leq h(t)$ for all $t\in\Rp$, from which
(\ref{eq:varphi}) follows.

The similar statement, involving $\varphi_-,\psi_-$ instead of $\varphi^-,\psi^-$, is obtained immediately by multiplying by $-1$.
\end{proof}

\vspace{0.3cm}
\noindent{\bf 6. Proof of (\ref{eq:what_to_show}).}

Now we turn to prove that $\left(f(t;x,y)\right)^- =
\left(f(t;x,y)\right)_-$, for all $t\in \Rp, x\in X,y\in X^*$.
That is, we turn to prove that $E = \Rp$. 

Consider the space $\cX = \Rp \times X \times X^*$ with
half-open$\times$norm$\times$norm topology. Recall that with the
half-open topology $\Rp$ is a Baire space. 
Denote the subspace $\cS \times X \times X^*$ by $\cX_0$. A
straightforward computation shows that $f$ is jointly continuous
on $\cX_0$. It then follows that $(t,x,y)\mapsto f(t;x,y)^-$ is
upper and $(t,x,y)\mapsto f(t;x,y)_-$ is lower semicontinuous on
$\cX$, which means that the sets
$$A_n := \{(t,x,y) \in \cX : f(t;x,y)^- - f(t;x,y)_- < 1/n \}$$
are all open and contain the dense set $\cX_0$. We conclude that
the set
$$A := \bigcap_{n=1}^\infty A_n = \{(t,x,y) \in \cX : f(t;x,y)^- - f(t;x,y)_- = 0 \}$$
is a dense $G_\delta$ in $\cX$.

By the results in
\cite{Kur}, section II.22.V (sometimes referred to as the
Kuratowski-Ulam Theorem. To apply this theorem we need the
separability assumption), it follows
that there is a dense $G_\delta$ (in the half-open topology) set $E' \subseteq \Rp$ of points $t$ for
which the set
$$A_t = \{(x,y) \in X \times X^* :  \left(f(t;x,y)\right)^- = \left(f(t;x,y)\right)_- \}$$
is residual, and, in particular, dense in $X \times X^*$. But by
Lemma \ref{lem:subspace}, $A_t$ is closed, so for all $t\in E'$,
$A_t = X \times X^*$. In other words, we obtain that $E$ contains a dense $G_\delta$ in $\Rp$ in the
half-open topology, and it follows that $E$ is residual in $\Rp$
in the standard topology (because every open set $U$ in the half
open topology contains an open set $V$ in the standard one, such
that $V$ is dense in $U$).

By the discussion following the definition of $E$, we can define
$T_t x$ for all $t\in E$ and all $x\in X$, consistently with the
definition of $T_t x$ for $t\in \cS$. For $s,t \in \cS$, we have
\bes f(t+s;x,y) = f(t;T_s(x),y). \ees It follows that for
$t\in\Rp$, $s\in\cS$, \be\label{eq:s+t1} \left(f(t+s;x,y)\right)^-
= \left(f(t;T_s (x),y)\right)^-, \ee and similarly for $f_-$. So
whenever $t\in E$ and $s\in\cS$, then $t+s$ is also in $E$. Now in
(\ref{eq:s+t1}) we put $\cS \ni s_n \searrow s \in E$, to get, for all $t\in \cS$,
\begin{align*}
\left(f(t+s;x,y)\right)^-
&= \lim_n \left(f(t+s_n;x,y)\right)^- \\
&= \lim_n f(t;T_{s_n} (x),y) \\
(*)&= y(T_t(T_s(x))) \\
&= f(t;T_s(x),y)%\\
%(**)&= \left(f(t;T_s(x),y)\right)^-
\end{align*} 
(equality $(*)$ follows from the fact that
$T_t$ is weakly continuous).
%and $(**)$ follows from the fact that
%$\cS\subseteq E$).
It follows using Lemma \ref{lem:varphipsi} that
$$\left(f(t+s;x,y)\right)^- \leq \left(f(t;T_s(x),y)\right)^- \,\, , \,\, s\in E, t \in \Rp .$$
Similarly,
$$\left(f(t+s;x,y)\right)_- \geq \left(f(t;T_s(x),y)\right)_- \,\, , \,\, s\in E, t \in \Rp .$$
In particular, if $s,t \in E$, then
$$ \left(f(t+s;x,z)\right)^- \leq \left(f(t;T_s(x),y)\right)^- = \left(f(t;T_s(x),y)\right)_- \leq \left(f(t+s;x,z)\right)_- .$$
Thus, $E$ is a semigroup.

But then $E$ must be $\Rp$. Indeed, for $0<r\in\Rp$, $r - E$ contains a dense $G_\delta$ in $[0,r]$, so it must
intersect $E$. Thus $r$ is a sum of two elements in $E$, and hence is in $E$. It follows that $E = \Rp$, and the proof is complete.
\end{proof}

\begin{remark}
\emph{Note that for Hilbert spaces the above result is trivial,
because weak continuity implies strong continuity at $0$:
$$\|T_t h - h \|^2 = \|T_t h\|^2 - 2 \Re \langle T_t h, h \rangle + \|h\|^2 \leq 2\|h\|^2 - 2 \Re\langle T_t h, h \rangle \rightarrow 0 $$ 
as $t\rightarrow 0$ (see, for example, \cite[Section I.6]{SzNF70}), and strong continuity at $0$ implies uniform strong
continuity
(this remark -- that is, the \emph{triviality} of the
result -- is not true, in our opinion at least, for nonlinear
semigroups).}
\end{remark}

One might ask where in the proof we used the reflexivity of $X$. Checking the proof, one can see that
we need both $X$ and $X^*$ to be separable (in order to use the Ulam-Kuratowski Theorem), and that we need $X$ to be
weakly sequentially complete. These two conditions turn out to be equivalent to having $X$ separable and reflexive.

Another condition one might question is the contractiveness of the semigroup. This condition is not essential, as the following result shows.
\begin{corollary}
Let $X$ and $\cS$ be as above, and let $T = \{T_s\}_{s\in \cS}$ be a semigroup of operators on $X$ such that
(\ref{eq:continuity}) holds. Then $T$ can be extended to a weakly continuous semigroup of operators over $\Rp$ if and only if there exist $M,a\geq0$ such that for all $t \in \cS$,
\be\label{eq:ineq}
\|T_t\| \leq Me^{at} .
\ee
\end{corollary}
\begin{remark}
\emph{Any semigroup bounded on all bounded subsets of $\cS$ will satisfy (\ref{eq:ineq}) for appropriate $M$ and $a$.
Assuming that each $T_s$ is weakly continuous, the above result also holds for nonlinear semigroups, with the extended semigroup being only \emph{right}-weakly continuous.}
\end{remark}
\begin{proof}
It is a well known result that any weakly continuous semigroup over $\Rp$ satisfies (\ref{eq:ineq}) for appropriate $M$ and $a$, and for all $t\in\Rp$. Thus, if $T$ can be extended to a semigroup over $\Rp$, it must satisfy (\ref{eq:ineq}).

Conversely, if $T$ satisfies (\ref{eq:ineq}), then one can define a new semigroup $U$ by
$$U_s = e^{-as}T_s \,\, , \,\, s\in\cS .$$
Now one defines a new norm on $X$ by
$$\|x\|_{\rm new} = \sup_{s\in \cS} \|U_s x\| ,$$
and with this norm $U$ is a contractive semigroup (this is a standard construction). (\ref{eq:continuity}) and (\ref{eq:ineq}) together imply that $\|\cdot \|_{\rm new}$ is equivalent to $\|\cdot\|$. One checks that the normed space $(X, \|\cdot\|_{\rm new})$ is a separable, reflexive Banach space. Thus, with this new norm, $U$ satisfies the assumptions of Theorem \ref{thm:weak}, so it can be extended. Then one puts
$$T_t = e^{at}U_t \,\, , \,\, t \in \Rp$$
to obtain the desired extension of $T$.
\end{proof}

\section{Closing remarks}
The main limitation of Theorem \ref{thm:weak} is the conditions imposed on $X$, which make
it inapplicable to other cases of interest. In particular, the motivation for this study was an
attempt to extend densely parametrized semigroups of (normal) completely positive maps on von Neumann algebras, 
and the theorem as it stands cannot be used in that setting. However, a solution to this problem (extension of 
densely parametrized semigroups on von Neumann algebras) appears implicitly in \cite[pages 37-38]{SeLegue} for the case of 
unit preserving semigroups of normal $*$-endomorphisms on $B(H)$, and for semigroups of positive normal linear maps that are not necessarily unit preserving it will appear in \cite{SH07}. In a recent paper \cite{L07b}, the first named author 
proved an analog of Theorem \ref{thm:weak}
for arbitrary Banach spaces (for the case of linear operators), 
using a different approach.

\section{Acknowledgments}
The authors wish to thank Daniel Reem, Simeon Reich and Baruch Solel for reading several preliminary versions 
of this manuscript and for offering their comments.

%%%%%%%%%%%%%%%%%%%%%%%%%%%%%%%%%%%%%%%%%%%%%%%%%%%%%%%%%%%%%%%%%%%%%%%%%%%%
% The bibliography
\bibliographystyle{amsplain}

\end{document}